\documentclass[11pt]{article}
\usepackage[nosumlimits]{amsmath}
\usepackage{amssymb,amsthm,url}

\newcommand{\pp}[2]{\frac{\partial #1}{\partial #2}} 
\newcommand{\dede}[2]{\frac{\delta #1}{\delta #2}}
\newcommand{\dd}[2]{\frac{\diff#1}{\diff#2}}

\newcommand{\del}{\partial}
\def\MM#1{\boldsymbol{#1}}
\DeclareMathOperator{\diff}{d\!}

\DeclareMathOperator{\curl}{curl}

\usepackage{amscd}
\usepackage{natbib}
\usepackage{helvet}
\usepackage{amsfonts}
\usepackage{algorithmicx}
\usepackage[noend]{algpseudocode}
\usepackage{algorithm}

\usepackage[margin=2cm]{geometry}
\usepackage{graphicx}
\usepackage{color}
\usepackage{caption}
\usepackage{subcaption}
\usepackage{authblk}

\begin{document}
\title{A vertical slice frontogenesis test case for compressible nonhydrostatic dynamical
cores of atmospheric models}

\author[1]{Hiroe Yamazaki}
\author[*,1]{Colin J. Cotter}

\affil[1]{Department of Mathematics, Imperial College London, UK}
\affil[*]{Correspondence to: \texttt{colin.cotter@imperial.ac.uk}}

\maketitle

\begin{abstract}
A new test case is presented for evaluating the compressible dynamical cores of the atmospheric models. The test case is based on a compressible vertical slice model that can be obtained by simple modification of a standard three dimensional compressible dynamical core. On the one hand, an advantage of the test case is that is quasi-2D, so it can be run quickly on a standard workstation, enabling rapid experimentation with numerical schemes and discretisation choices. On the other hand, the test case exhibits frontogenesis, a challenging regime for numerical discretisations which usually only arises in 3D model configurations for the compressible case. Numerical results of the test case using an implicit time-stepping method with a compatible finite element discretisation are presented as a reference solution. An example comparison between advective and vector-invariant forms for the advective nonlinearity in the velocity equation demonstrates one possible use of the scheme. The comparison shows a Hollingsworth-like instability when the vector invariant form is used.
\end{abstract}

\noindent \textbf{keywords:} frontogenesis; vertical slice; numerical weather prediction; model intercomparison

\section{Introduction}

Standardised test cases on hierarchies of models are a crucial element of atmospheric model dynamical core development. They are useful for code verification (checking that the code is a correct implementation of the mathematics), evaluation of the quality of discretisation choices (considering e.g. conservation properties, numerical dissipation and dispersion effects, etc.), benchmarking of computational speed and parallel scalability of models, and general intercomparison of models between modelling groups \citep{williamson1992standard,galewsky2004initial,lauter2005unsteady,ullrich2015analytical,ullrich2012dynamical, hall2016dynamical, ullrich2017dcmip2016, zarzycki2019dcmip2016}.

Ideally, we would like to use testcases to develop confidence that the numerical solutions are converging to the true solution, but error calculations require a known exact solution. Exact solutions of the nonlinear equations are hard to find, especially solutions that are relevant to specific types of atmospheric flow. As an alternative, \citet{cullen2007modelling} proposed to consider asymptotic limits (such as the quasigeostrophic or semigeostrophic limit), demonstrating convergence in the asymptotic limit to numerical solutions of the limit equation where the error is more under control. 

Here we focus on the vertical slice Eady problem, which is important because it demonstrates the ability to simulate the geostrophic flows which are the basis of extratropical weather. In the incompressible Boussinesq case, the Euler solutions can be shown formally to converge to the solutions of the simpler semi-geostrophic system as a scaling parameter tends to zero. The latter solutions are discontinuous, and can be rigorously computed to any desired accuracy by optimal transport methods \citep{cullen2006mathematical,egan2022new,benamou2024entropic, lavier2024semi}. They thus provide a challenging test for numerical methods for the Euler equations. Further, considering Euler, any solution will require regularisation, and \citet{nakamura1994nonlinear} showed strong sensitivity to the choice of regularisation. The semigeostrophic limit solution is well defined, and conserves energy, with no need for regularisation. Hence, it can be argued that the correct limiting solution for Euler should have minimal numerical regularisation (explicit or implicit). For the incompressible vertical slice Eady model in the semigeostrophic limit, the asymptotic limit programme was pursued in \citet{visram2014framework, yamazaki2017vertical}.

The incompressible Boussinesq equations are not an ideal test for atmospheric models as they are not a good approximation to the compressible Euler equations, which are used in production atmospheric models, in problems with a significant vertical extent. There is thus a strong desire to reformulate the Eady test case for the compressible equations.

 In this paper we present a compressible vertical slice model with initial conditions and parameter settings that lead to a frontogenesis lifecycle similar to the solutions considered in the incompressible case, and we offer this to the modelling community as an additional tool in their evaluation and intercomparison of atmosphere models.

The paper is structured as follows. Section 2 presents the full details of the new frontogenesis test case of compressible flow, including formation
of the compressible Eady problem, settings of the frontogenesis experiment and the initial conditions.
In Section \ref{results}, we apply the test case to an example model using an implicit time-stepping method with a
compatible finite element discretisation, in order to provide a reference solution for other models of different numerical schemes. To evaluate the validity and reliability of the new test case, the results are compared to the solutions from the incompressible Euler–Boussinesq Eady models.
A summary and outlook is provided in Section \ref{conclusion}.

\section{Test case description}\label{model_description}
\subsection{Governing equations}

In this section, we describe the model equations that are employed for the new vertical slice test case for compressible models. In the incompressible case, a 3D Eady problem can be solved using 2D computation by assuming that all variables are independent of $y$, except for temperature which has a time independent linear dependency on $y$, necessary to support the baroclinic instability that leads to frontogenesis. When compressibility is introduced, it is no longer possible to find similar quasi-2D solutions of the 3D model. We can only introduce \emph{ad hoc} terms (motivated by the effects of a $y$-direction temperature gradient) to produce baroclinic instability. This was attempted in \cite{cullen2008comparison}, but unfortunately the proposed modifications did not have a conserved energy or potential vorticity, which are crucial in sustaining baroclinic instability towards frontogenesis. The question of how to build vertical slice models with conserved energy and potential vorticity was addressed by \citet{cotter2013variational}, who identified the variational structure underpinning the incompressible vertical slice model. This structure establishes energy and potential vorticity conservation from within the vertical slice setting without needing to refer to the underlying 3D model. Then, other models that also conserve energy and potential vorticity can be derived by choosing a different Lagrangian (kinetic energy minus potential energy). \citet{cotter2019particle} showed that the potential vorticity conservation in these models arises from a special form of particle relabelling symmetry in the slice geometry that incorporates the velocity
in the $y$-direction. \citet{cotter2013variational} used the variational structure to propose other vertical slice Eady models, including a compressible model. This zoo was extended to anelastic and pseudocompressible Eady models in \citet{cotter2014variational}. When moving away from a linear equation of state, models in this framework are not a reduction of a 3D counterpart, and they serve only to emulate an Eady flow configuration in the vertical slice setting for the purpose of evaluating numerical models. In the present work, we investigated the compressible vertical slice Eady model from this framework. It was hoped that the energy and potential vorticity conservation would enable sustained baroclinic instability until frontogenesis, and we have indeed found this to be the case. We have made one minor modification to the model proposed by \citet{cotter2013variational}, adding a constant perturbation $\Pi_0$ to the Exner pressure to reduce the vertical mean horizontal flow, making the setup more suitable for a numerical test problem.

The equations are
\begin{eqnarray}
\frac{\partial \MM{u}}{\partial t} + (\MM{u} \cdot \nabla)\MM{u} - f v \, \MM{\hat{x}} &=& -c_p \theta_S \nabla \Pi - g \MM{\hat{z}}, \label{ueq}\\
\frac{\partial v}{\partial t} + (\MM{u} \cdot \nabla){v} + f u &=& \underbrace{c_p s(\Pi-\Pi_0)}_{\mbox{extra}}, \label{veq}\\
\frac{\partial \theta_S}{\partial t} + (\MM{u} \cdot \nabla)\theta_S &=& \underbrace{- vs}_{\mbox{extra}}, \label{thetaeq}\\
\frac{\partial D}{\partial t} + \nabla \cdot (D \MM{u}) &=& 0, \label{rhoeq}
\end{eqnarray}
where $\MM{u} = (u, w)$ is the velocity vector in the vertical slice, $\nabla = (\partial_x, \partial_z)$ is the gradient operator in the vertical slice, $v$ is the out-of-slice velocity component, $\MM{\hat{x}}$ is a unit vector in the $x$-direction, and $\MM{\hat{z}}$ is that in the $z$-direction; $D$ is the density, $\theta$ is the potential temperature, $f$ is the rotation frequency, $g$ is the acceleration due to gravity, and $\theta_S$ is a ``perturbation potential temperature",
with the idea that the total potential temperature $\theta\approx sy + \theta_S$ in the region near
$y=0$ where the slice model is defined, and for scaling purposes we write
\begin{eqnarray}
s &=& - \frac{\theta_0 f \Lambda}{g} = \text{const.},
\end{eqnarray}
where $\theta_0$ is a reference density at the surface, and $\Lambda$ is the constant vertical shear.
Further, the Exner pressure $\Pi$ depends on $\theta_S$ and $D$ such that
\begin{equation}
\Pi = \left(\frac{p}{p_0}\right)^{R/c_p}, \,
p = DR\underbrace{T}_{=\theta_S\Pi},
\label{eq:exner}
\end{equation}
where $p$ is the pressure, $p_0$ is a reference pressure, $R=c_p-c_v$ is the specific gas constant, $c_p$ is the specific heat at constant pressure, and $c_v$ is the specific heat at constant volume.
All the variables $\MM{u}$, $v$, $\theta_S$, $D$ and $\Pi$ depend on $(x, z, t)$ only, hence the model
can be solved on a 2D mesh.

It is important to note here that it is not physically corrrect to have $\Pi$ depending on $\theta_S$
rather than $\theta$, but this modification is necessary to have a well-defined vertical slice model. This model is not intended to skillfully predict physical phenomena, but rather to provide a challenging
test case for atmosphere model numerics in a 2D configuration. This configuration can also be obtained
by setting periodic boundary conditions to a 3D model (in local area model mesh configuration) in the $y$-direction and then making the mesh only 1 element wide (which is what we did in this paper). This is then testing the same code that is used in the 3D model, which is useful for verification purposes; it is also potentially more achievable in a code framework that does not easily switch between 3D and 2D. After this is done, the only other required changes to a 3D model are to insert the terms indicated as ``extra" in Equations (\ref{ueq}-\ref{rhoeq}), having identified $\theta$ in the 3D model with $\theta_S$ in the vertical slice model. 

The vertical slice model has steady state solutions when $\theta_S$ depends on $z$ only, with $v=0$,
$\Pi$ obtained from the hydrostatic balance,
\begin{equation}
0 = -c_p\theta_S\pp{\Pi}{z} - g,
\end{equation}
and $u$ obtained from the geostrophic balance, 
\begin{equation}
fu = c_ps(\Pi-\Pi_0).
\end{equation}
This solution can only be sustained with periodic boundary conditions in the $x$-direction.
Here we see the role of $\Pi_0$: without it, the vertical mean of $u$ will be large and everything is propagating rapidly to the right.

This modification has a conserved total energy $E = K_{u} + K_{v} + P$, given by
\begin{eqnarray}
K_{u} &=& \frac{1}{2} \int_{\Omega} D |\MM{u}|^{2} \, \diff x \diff z, \label{kinetic_u}\\
K_{v} &=& \frac{1}{2} \int_{\Omega} D v^{2} \, \diff x \diff z,\label{kinetic_v}\\
P &=& \int_{\Omega} D \, (gz + c_v \Pi \theta - c_p \Pi_0 \theta) \, \diff x\diff z, \label{potential}
\end{eqnarray}
where $K_u$ is the kinetic energy from the in-slice velocity components, $K_v$ is the kinetic energy from the 
out-of-slice velocity component, and $P$ is the potential energy comprising contributions from 
gravitational potential energy and internal energy. For further discussion of the energy
conservation and potential vorticity, see the Appendix.

\subsection{Experimental settings}
\subsubsection{Constants}\label{constants}

In the frontogenesis experiments, the model constants are set to the values below, 
following the incompressible
Euler–Boussinesq Eady test case \citep[e.g.][]{nakamura1994nonlinear, 
cullen2007modelling, visram2014framework, visram2014asymptotic, yamazaki2017vertical},
\begin{eqnarray}
L &=& 1000\ \mathrm{km},\ H = 10\ \mathrm{km},\ f = 10^{-4}\ \mathrm{s}^{-1},
\ \nonumber\\
g &=& 10\ \mathrm{m\ s}^{-2},\ p_{0} = 10^5 \ \mathrm{Pa},\ \theta_{0} 
= 300\ \mathrm{K}, \nonumber\\
\Lambda &=& 10^{-3}\ \mathrm{s}^{-1},\ N^{2} = 2.5 \times 10^{-5}\ \mathrm{s}^{-1},
 \nonumber \\
 \Pi_0& =& 0.864,
\end{eqnarray}
where $L$ and $H$ determine the model domain $\Omega = [-L, L] \times [0,H]$, and $N$ is the Brunt-V\"ais\"al\"a frequency. 
The $y$ derivative, $s$, of the background temperature in \eqref{thetaeq} and \eqref{veq} 
is therefore calculated as
\begin{eqnarray}
s = - \frac{\theta_0 f \Lambda}{g} = - 3.0 \times 10^{-6}. 
\end{eqnarray}
The Rossby and Froude numbers are given in the model as
\begin{eqnarray}
\mathrm{Ro} = \frac{u_{0}}{fL} = 0.05, \label{Rossby} \\
\mathrm{Fr} = \frac{u_{0}}{NH} = 0.1,
\end{eqnarray}
where $u_{0} = 5\ \mathrm{m\ s}^{-1}$ is used as a representative velocity. 
The ratio of Rossby number to the Froude number defines the Burger number,
\begin{eqnarray}
\mathrm{Bu} = \mathrm{Ro}/\mathrm{Fr} = 0.5, 
\end{eqnarray}
which is used in the next section for initialising the model.

\subsubsection{Initialisation}\label{initialisation}

The model field is initialised with a small perturbation to a balanced steady state.
The steady state temperature is isothermal,
\begin{equation}
\bar{\theta}(z) = \theta_{0}\exp(N^2(z-H/2)/g).
\end{equation}
Then,
following the incompressible
Euler–Boussinesq Eady test case \citep[e.g.][]{nakamura1994nonlinear, 
cullen2007modelling, visram2014framework, yamazaki2017vertical}, a small perturbation is applied to the in-slice temperature, 
\begin{eqnarray}
\label{perturbation}
\theta_S(x, z) &=& \bar{\theta}+\frac{\theta_0 a N}{g} \left\{-\left[1-\frac{\mathrm{Bu}}{2}\coth\left(\frac{\mathrm{Bu}}{2}\right)\right] 
\sinh Z \cos \left( \frac{\pi x}{L} \right) -n \,\mathrm{Bu} \cosh Z \sin 
\left( \frac{\pi x}{L} \right) \right\}, 
\label{init_buoyancy}
\end{eqnarray}
which is the structure of the normal mode taken from \citet{williams1967atmospheric}. 
The constant $a$ corresponds to the amplitude of the perturbation,
and the constant $n$ takes the form of
\begin{eqnarray}
n = \frac{1}{\mathrm{Bu}}\left\{ \left[\frac{\mathrm{Bu}}{2} - 
\tanh \left(\frac{\mathrm{Bu}}{2} \right) \right] 
\left[\coth \left( \frac{\mathrm{Bu}}{2} \right) - \frac{\mathrm{Bu}}{2}  \right] \right\}^{\frac{1}{2}}.
\end{eqnarray}
The modified vertical coordinate $Z$ is defined as
 \begin{eqnarray}
Z = \mathrm{Bu}\left[\left(\frac{z}{H} \right) - \frac{1}{2} \right].
\end{eqnarray}
In fact, this is not a normal mode of the linearisation of the compressible vertical slice model.
Instead, we rely upon initialising a balanced initial condition which will trigger a balanced
baroclinically  unstable mode. This is why our experiment  has an initial growth period before
resetting the clock to 0 when the amplitude starts to become visible.

Next we initialise the density $D$. Given the initial temperature $\theta_S$ in \eqref{perturbation}, 
first we seek $\Pi$ in hydrostatic balance,
\begin{equation}
\label{pvbalance}
c_p \theta_S \frac{\partial \Pi}{\partial z} = -g.
\end{equation}
Once we have the profile of the initial  $\Pi$, we can calculate the initial density by using \eqref{eq:exner} as
\begin{equation}
\label{pvbalance2}
D = \frac{p_0 \Pi^{\frac{c_v}{R_d}}}{R_d \theta_0}.
\end{equation}
Alternatively, we can solve the nonlinear problem of finding $D$ such that \eqref{pvbalance}
holds, treating $\Pi$ as the prescribed function of $D$ and $\theta_S$ given in 
\eqref{eq:exner}. This is what we did in our experiments (see \citet{cotter2023compatible} for details).

Then we initialize the out-of-slice velocity $v$ by seeking a velocity in geostrophic balance with the initialised $\theta_S$ and $\Pi$:
\begin{eqnarray}
- f v = - c_p \theta_S \frac{\del \Pi}{\del x}. 
\end{eqnarray}

Finally, we initialise the in-slice velocity ${u}$ according to
\begin{eqnarray}
u &=& \frac{c_p}{f}\frac{\del \bar{\theta}}{\del y} (\Pi-\Pi_0).
\end{eqnarray}
Here $\Pi_0$ was chosen following numerical experiments to minimise the horizontal drift
of the front structures once they have formed.

In addition to the initialisation of the model
field described above, we introduce a breeding procedure at the beginning of each simulation to remove any remaining unbalanced modes in the initial condition, as in \citet{visram2014framework, visram2014asymptotic, yamazaki2017vertical}.
In the experiments performed in section \ref{results}, 
the model field is initialised with a small perturbation by choosing $a$ = -7.5 in 
\eqref{init_buoyancy}, following the previous studies. The simulation is then advanced 
until the maximum amplitude of $v$ reaches 3 m s$^{-1}$, 
at which point the time is reset to zero to match the amplitude of the initial perturbation 
with that of  the incompressible Euler–Boussinesq Eady models \citep[][]{nakamura1989nonlinear, visram2014framework, yamazaki2017vertical}
as closely as possible. In the results presented in Section \ref{results}, the time is reset to zero after 50 hours, both with the default resolution and with the high resolution.  

\section{Results}\label{results}
In this section, we present the results of the frontogenesis experiments using the 
new test case developed in this study, with the use of the finite element 
code generation library Firedrake \citep{FiredrakeUserManual}. Since this paper
is about the testcase and not the specific numerical discretisation, we just 
briefly mention that we used a compatible finite element discretisation following
\citet{cotter2023compatible} in
3D configuration on a mesh that is one cell wide and periodic in the $y$-direction, with just two differences. First, we did not 
use the edge stabilisation recommended there for $\theta_S$ transport as it led
to instabilities after the front formed, requiring further investigation. It appears
that the horizontal discontinuous Galerkin upwinding for $\theta_S$ is sufficient
for this problem, probably because the curvature of the solution is mainly in the
horizontal direction. Second, we did not use a vector invariant form of the velocity
advection equation, instead using a discretisation of the standard advective form
$(\MM{u}\cdot\nabla)\MM{u}$
in our main results, using the upwind formulation of \citet{cockburn2007note} for
that term instead. This is discussed more later. The timestepping
method is the implicit midpoint rule, which is solved using Newton's method and 
GMRES preconditioned by the column patch additive Schwarz method, which is also 
described in \citet{cotter2023compatible}. This scheme is unconditionally stable
and preserves energy for linear wave equations; any dissipation in the model is 
coming from upwinding in the spatial discretisation.

The constants used to set up the experiments are shown in section \ref{constants}. 
At the beginning of each experiment, the model is initialised in the way described 
in section \ref{initialisation},
then integrated for 25 days in each experiment. 

In our solution, the model resolution is given by
\begin{eqnarray}
\Delta x = \frac{2L}{N_{x}},\ \Delta z = \frac{H}{N_{z}}, \label{grid_space}
\end{eqnarray} 
where $N_x$ and $N_z$ are the number of quadrilateral elements in the $x$- 
and $z$-directions, respectively. 
Unless stated otherwise, we use a resolution of $N_{x}$ = 30 and $N_{z}$ = 30, and the time step $\Delta t$ of 300 s.

\begin{figure}[p]
  \centering
  \begin{subfigure}{0.45\hsize}
    \centering
    \includegraphics[width=85mm]{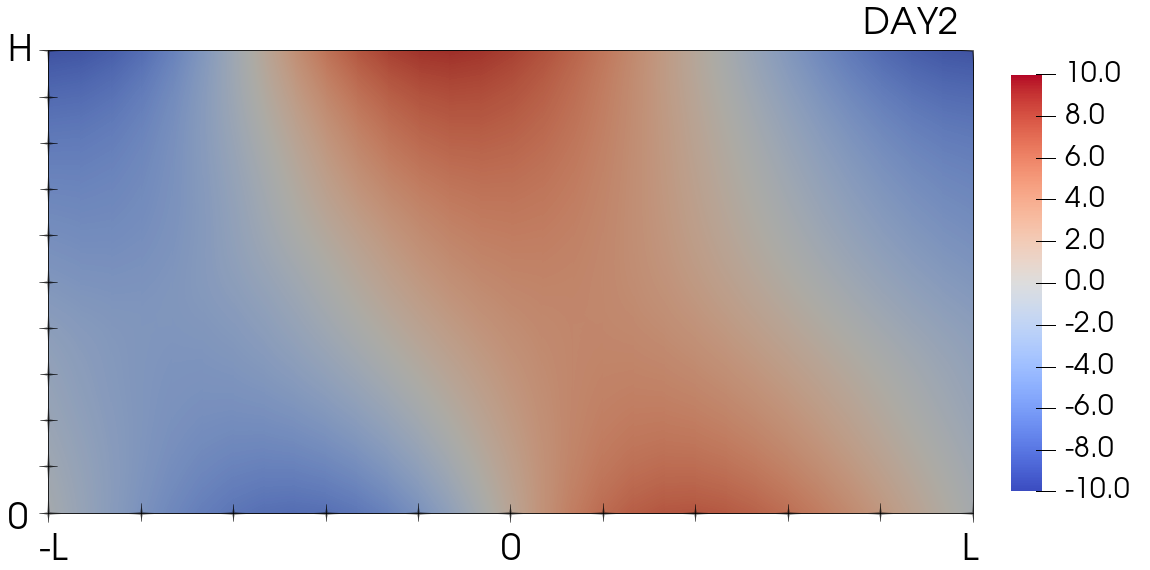}
    \includegraphics[width=85mm]{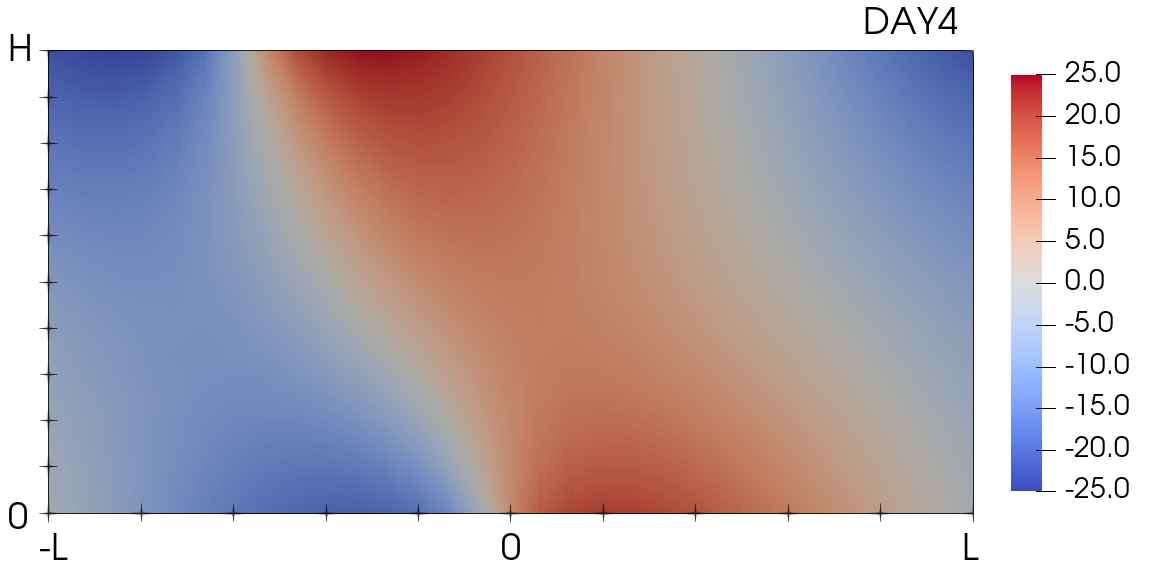}
    \includegraphics[width=85mm]{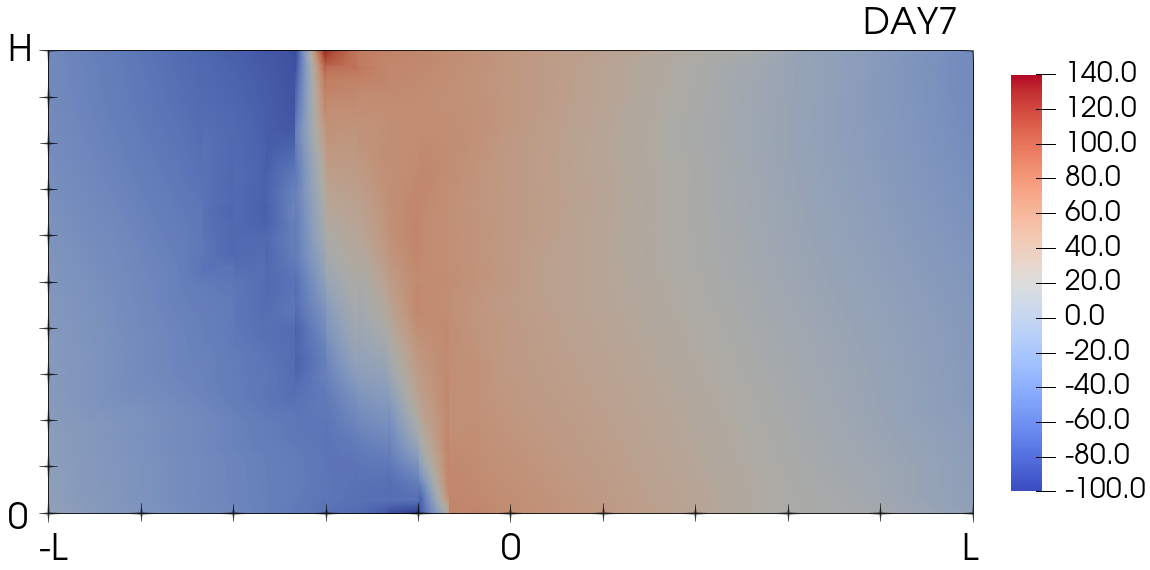}
    \includegraphics[width=85mm]{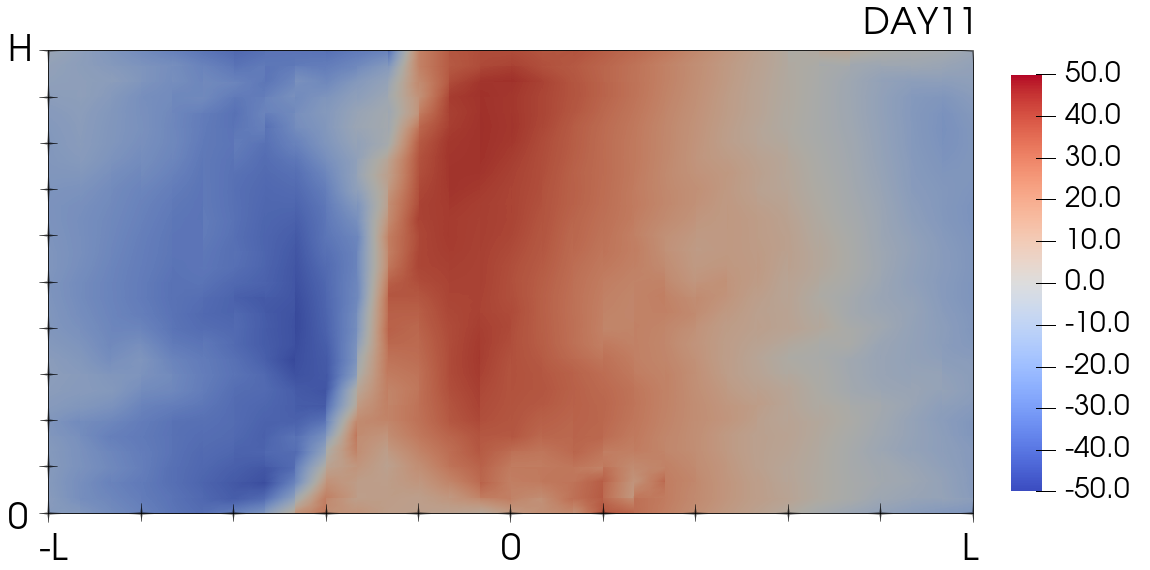}
    \caption{control run}
    \label{fig:velocity-contours-control}
  \end{subfigure}
  \hspace{2em}
  \begin{subfigure}{0.45\hsize}
    \centering
    \includegraphics[width=85mm]{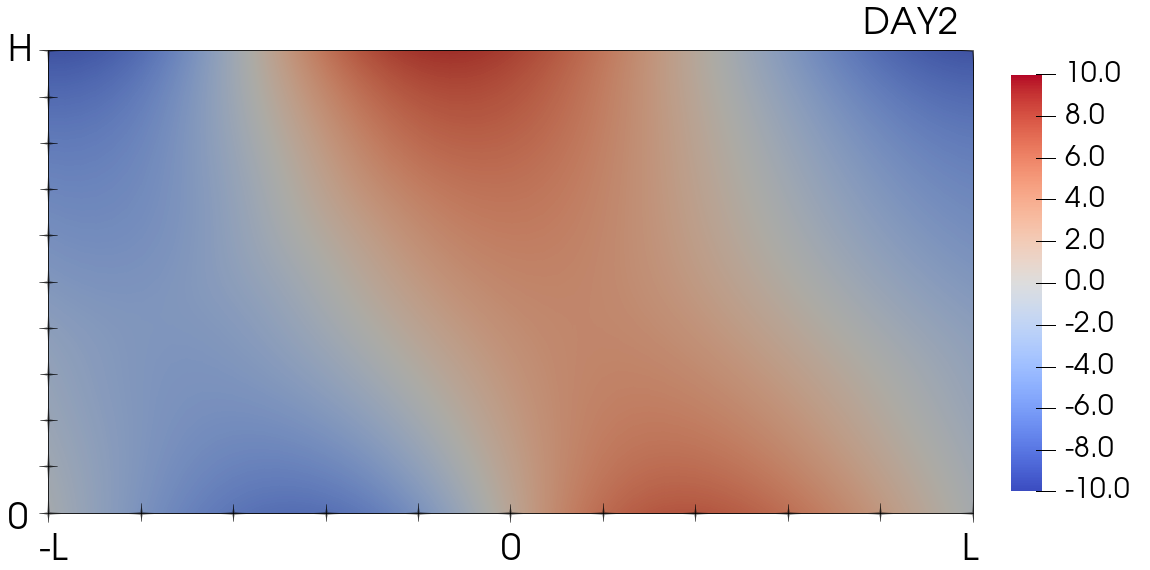}
    \includegraphics[width=85mm]{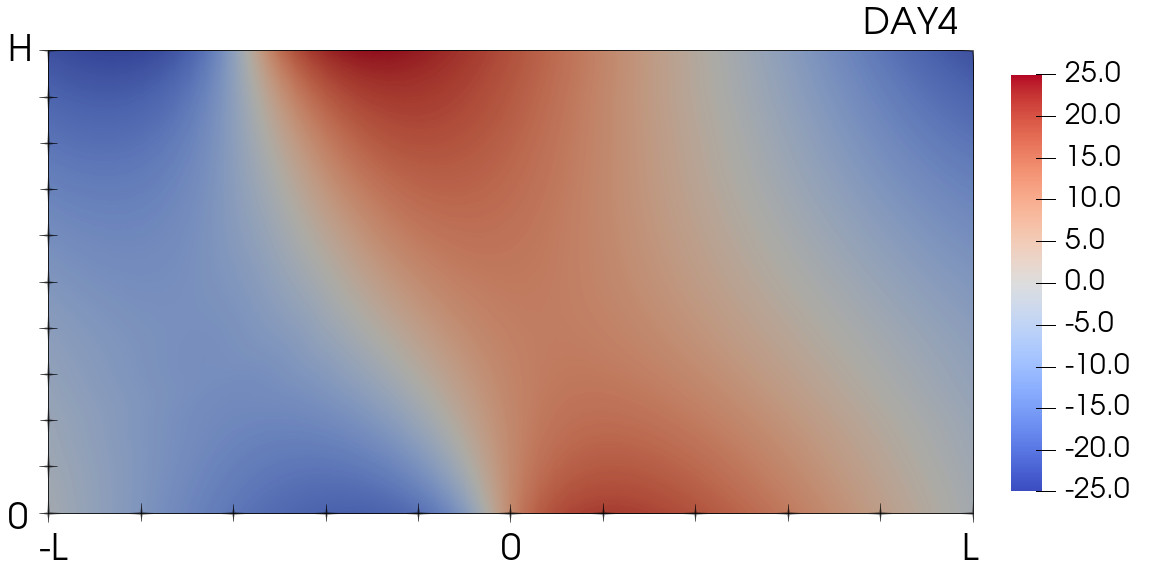}
    \includegraphics[width=85mm]{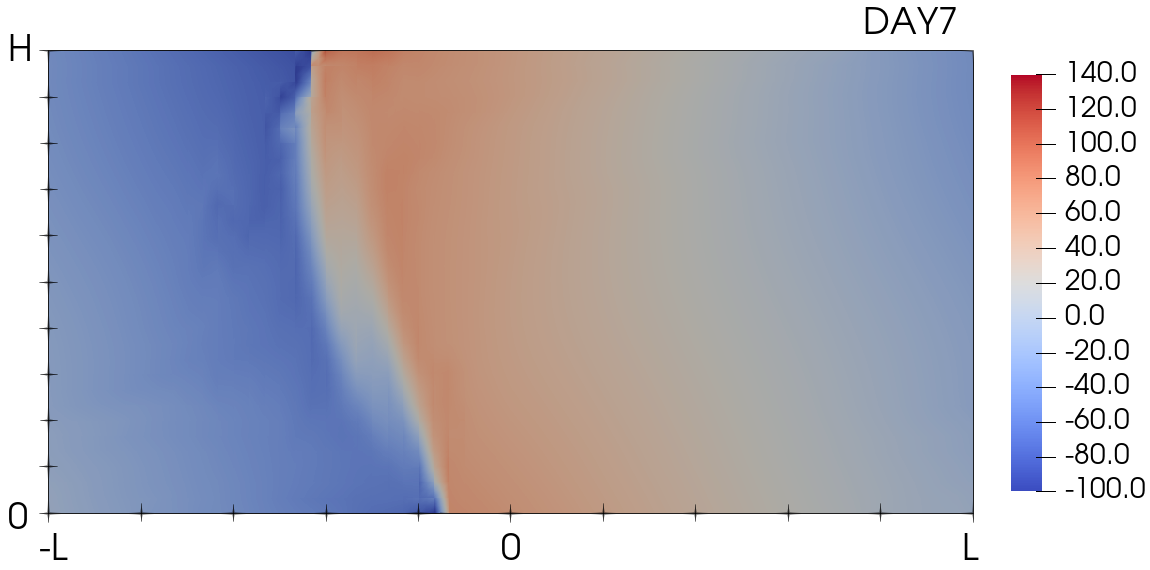}
    \includegraphics[width=85mm]{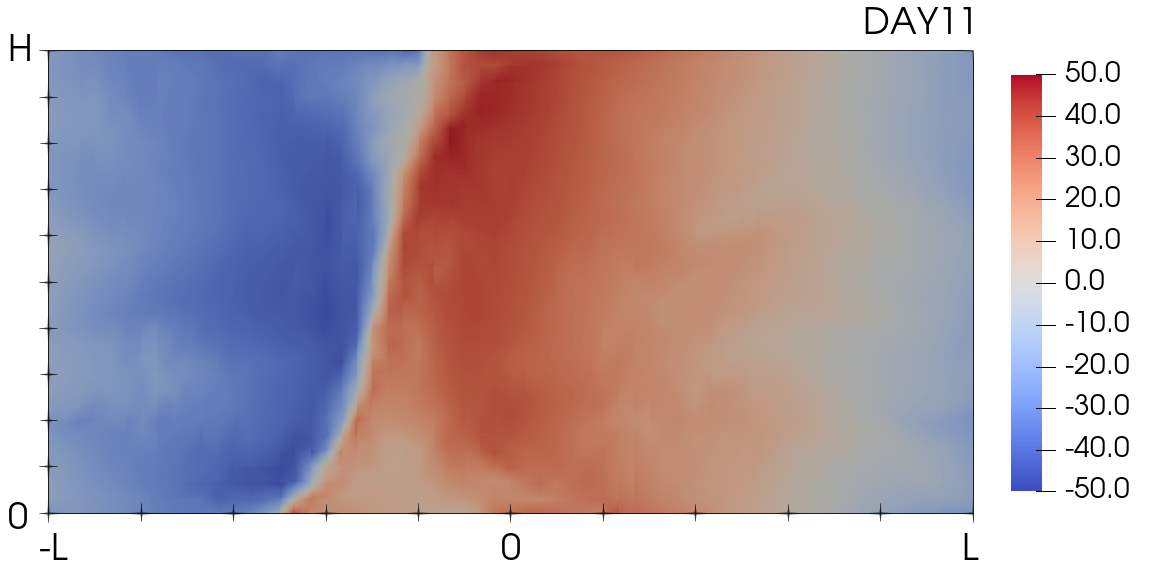}
    \caption{high-resolution run}
    \label{fig:buoyancy-contours-highres}
  \end{subfigure}
  \caption{Snapshots of out-of-slice velocity fields at days 2, 4, 7, and 11 in (a) the control run and (b) the high-resolution run.}
  \label{1fig:slice-contours v}
\end{figure}

\begin{figure}[p]
  \centering
  \begin{subfigure}{0.45\hsize}
    \centering
    \includegraphics[width=85mm]{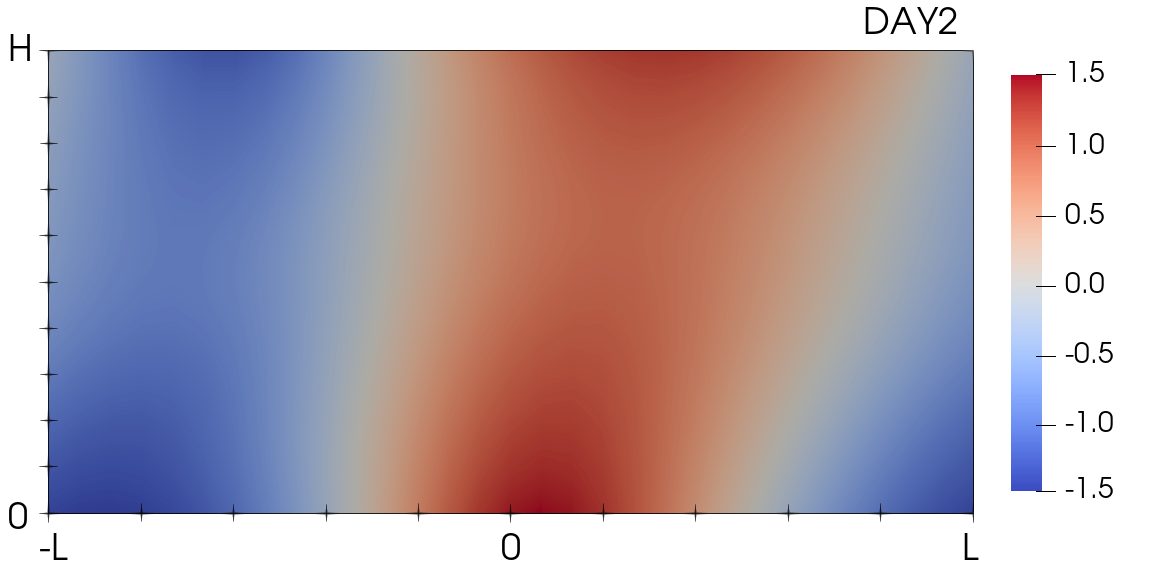}
    \includegraphics[width=85mm]{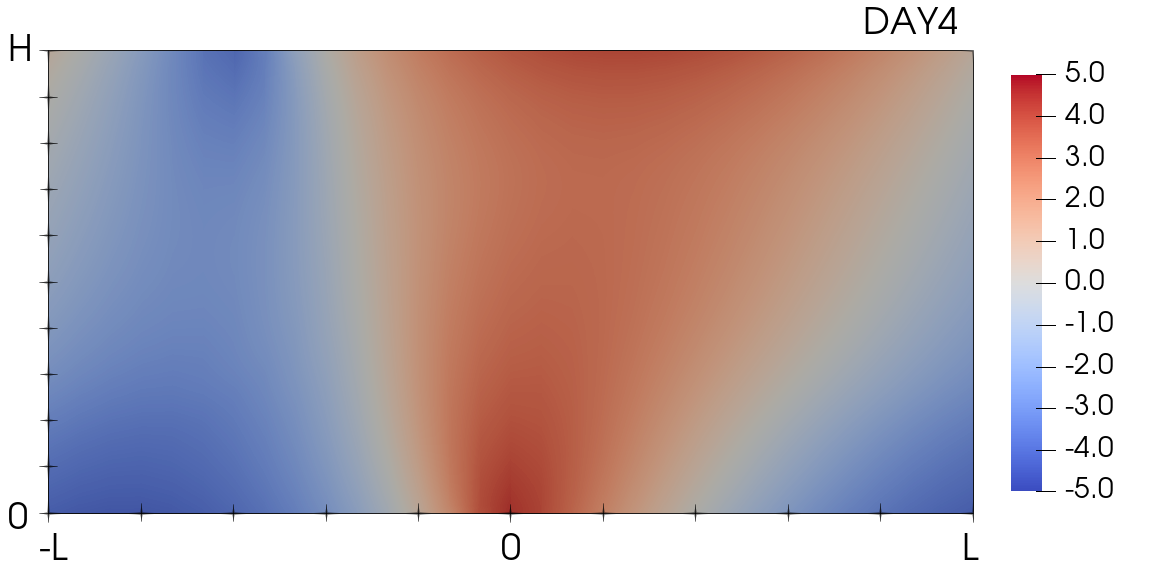}
    \includegraphics[width=85mm]{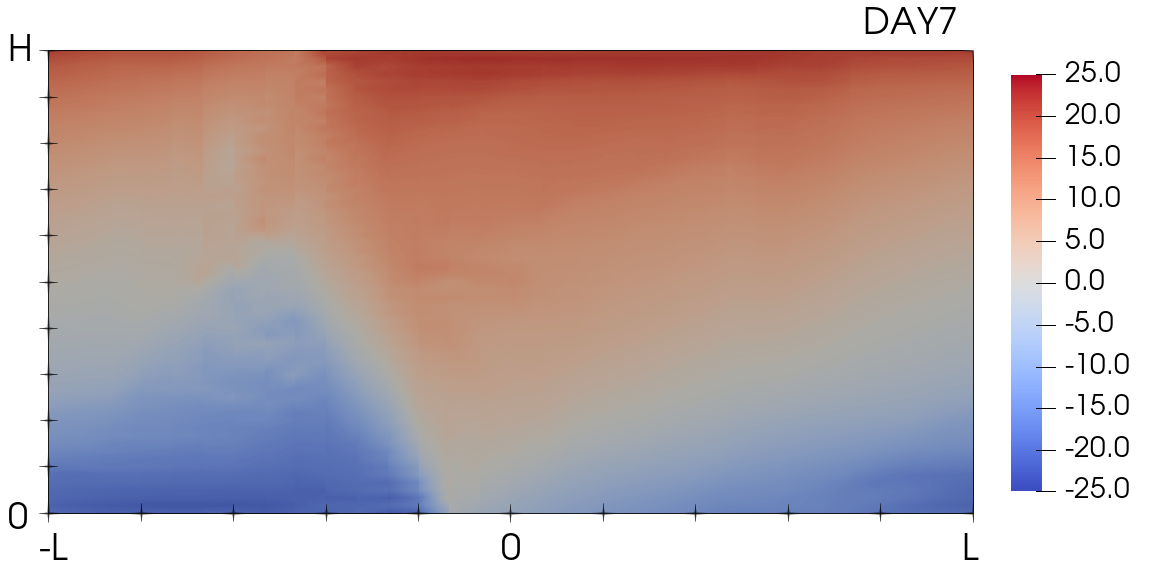}
    \includegraphics[width=85mm]{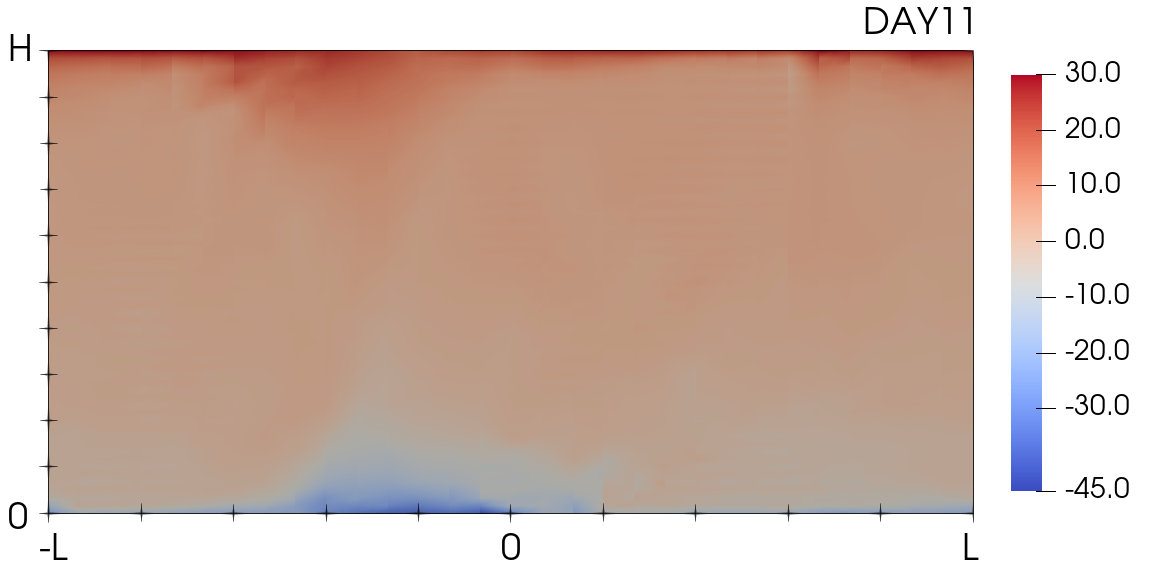}
    \caption{control run}
    \label{fig:theta-contours-c}
  \end{subfigure}
  \hspace{2em}
  \begin{subfigure}{0.45\hsize}
    \centering
    \includegraphics[width=85mm]{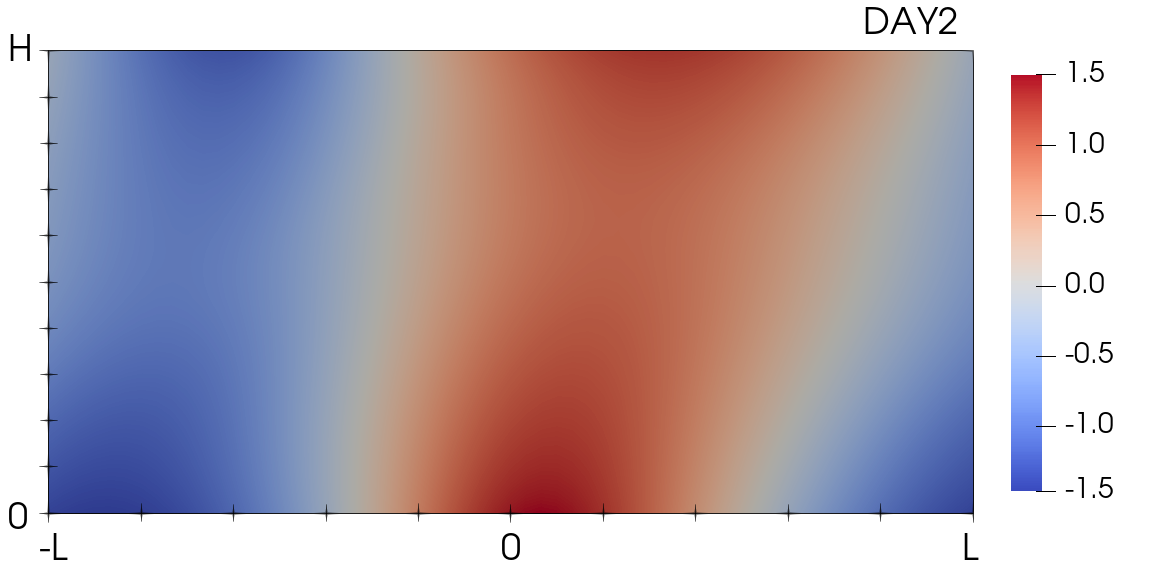}
    \includegraphics[width=85mm]{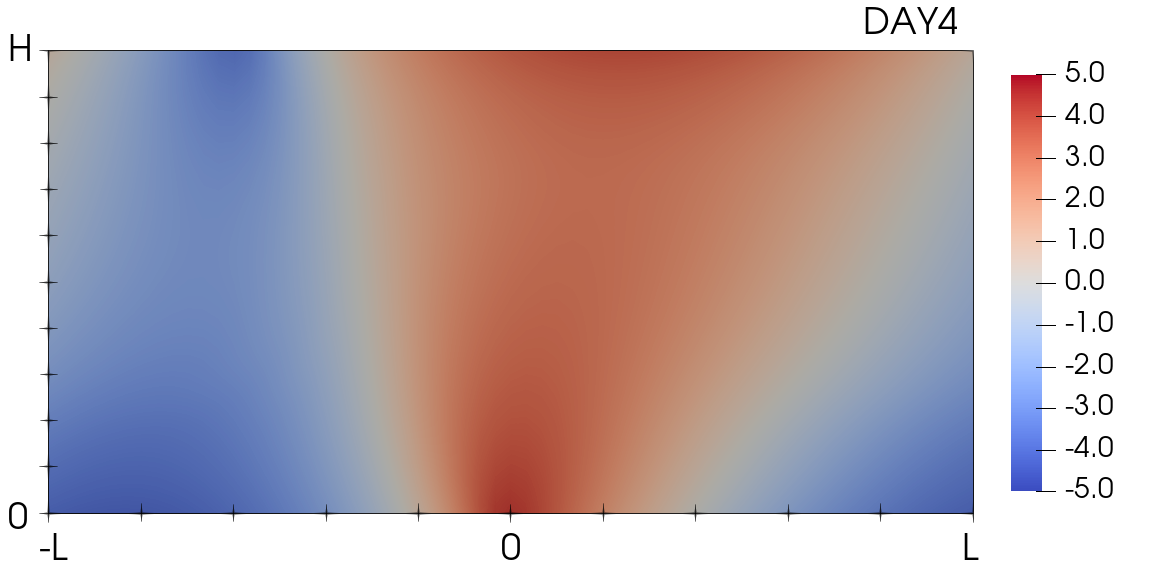}
    \includegraphics[width=85mm]{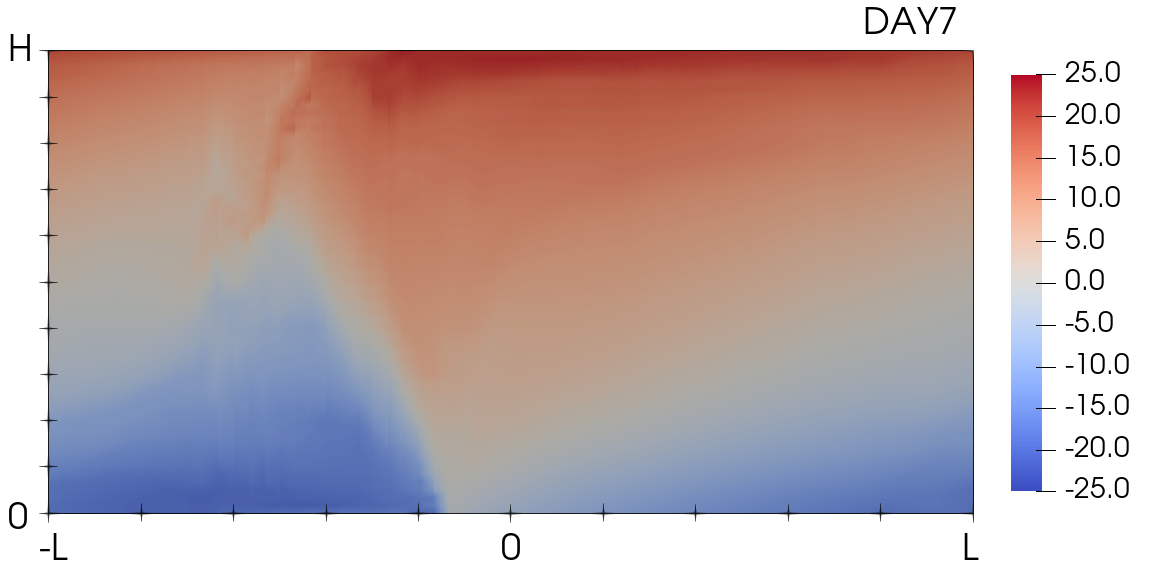}
    \includegraphics[width=85mm]{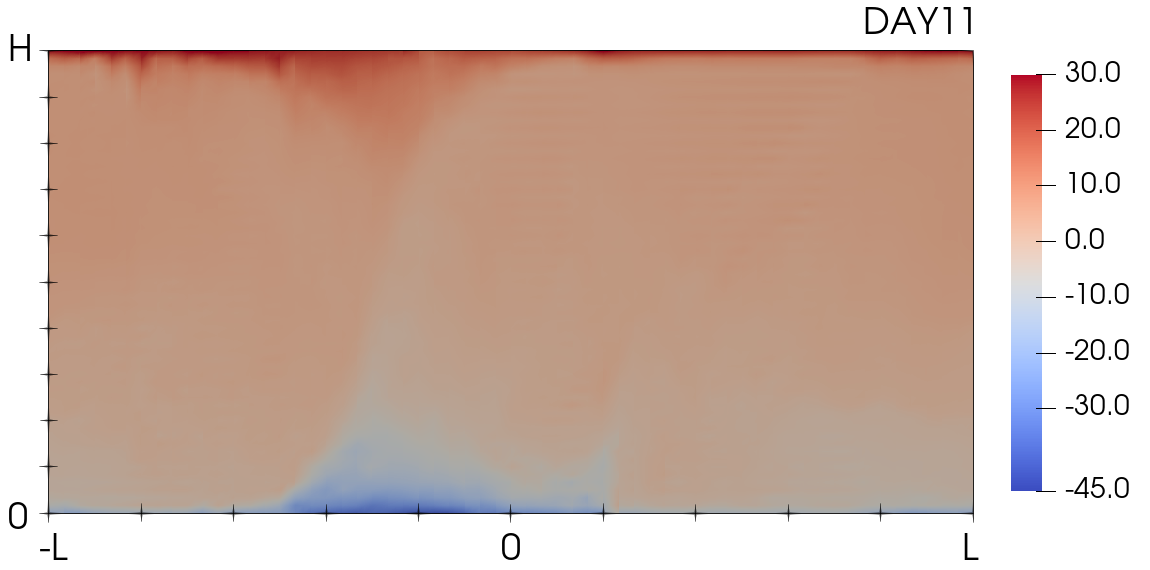}
    \caption{high-resolution run}
    \label{fig:buoyancy-contours-h}
  \end{subfigure}
  \caption{Snapshots of in-slice temperature perturbation field at days 2, 4, 7, and 11 in (a) the control run and (b) the high-resolution run.}
  \label{1fig:slice-contours theta}
\end{figure}

\begin{figure}[htbp]
  \centering
  \includegraphics[width=120mm,angle=0]{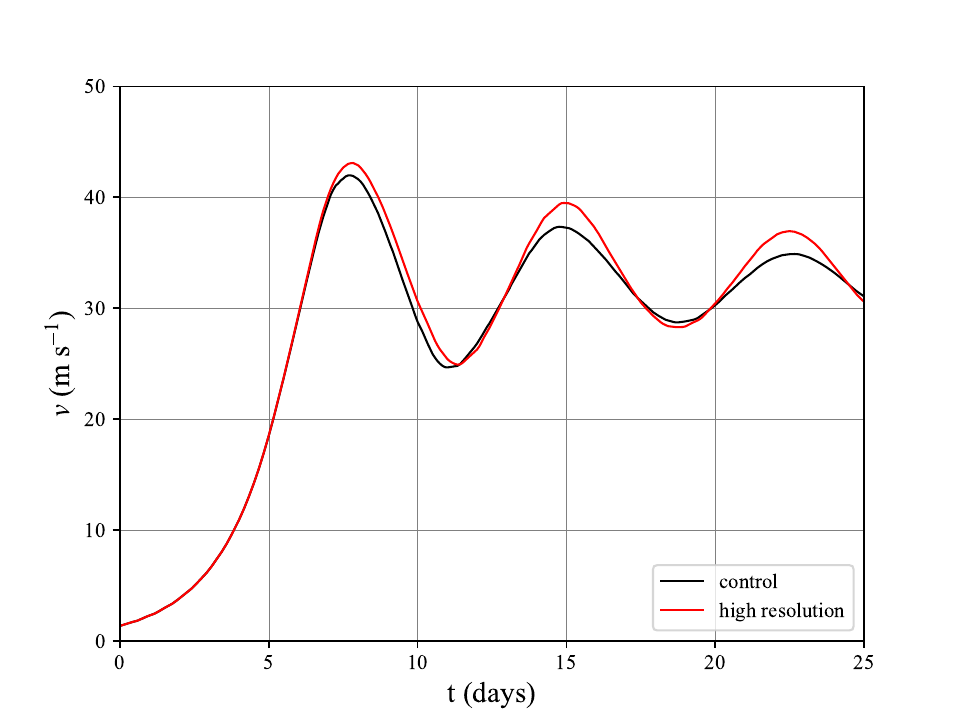}
  \caption{The root mean square of the out-of-slice velocity in the control run (black) and in the high-resolution run (red).}
  \label{fig:rmsv-velocity-comparison}
\end{figure}

Figures \ref{1fig:slice-contours v} and \ref{1fig:slice-contours theta} show the snapshots 
of the out-of-slice velocity $v$ and potential temperature $\theta$ fields, respectively. 
At day 2, both fields show very similar structures to those from the simulation using 
the linearised equations in \citet{visram2014asymptotic}. 
The model shows some early signs of front formation at day 4. 
The frontal discontinuity becomes most intense around day 7. 
At day 11, the vertical tilt in the $v$-field reverses, which is a sign of energy 
conversion from kinetic back to potential energy.
Overall, these results are consistent with the incompressible Eady models 
\citep[][]{nakamura1989nonlinear,
cullen2007modelling, visram2014framework, visram2014asymptotic, yamazaki2017vertical}. 
The black curve in Figure \ref{fig:rmsv-velocity-comparison} shows the time 
evolution of the root mean square of $v$ (RMSV) in our model. 
The result shows that the model reproduces several further quasi-periodic lifecycles 
after the first frontogenesis. 

To evaluate the effect of resolution on our model result, 
we repeated the experiment using $N_x$ = 60 in stead of $N_x$ = 30. A time step of $\Delta t$ = 120 s is used for the high-resolution run. The number of layers is kept at $N_z$ = 30.
The evolution of RMSV in the high-resolution run is shown by the red curve in Figure 
\ref{fig:rmsv-velocity-comparison}. It shows front cycles with bigger RMSV at the peaks compared the the low-resolution results. 

\begin{figure}[htbp]
  \centering
  \includegraphics[width=120mm]{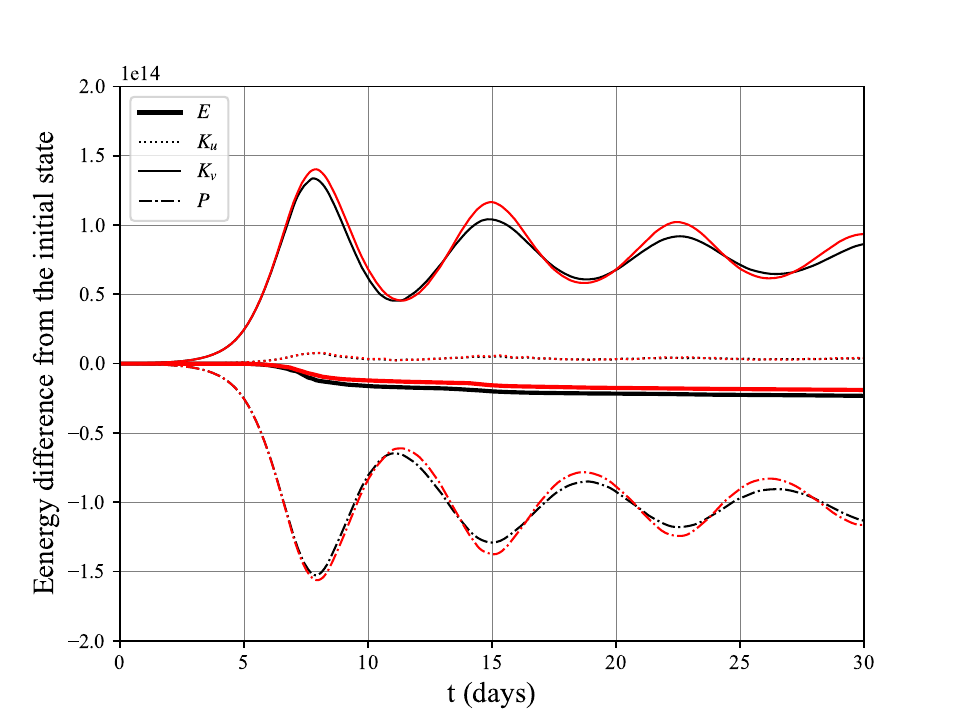}
  \caption{Time evolution of energy in the control run (black) and in the high-resolution run (red).
  Thick line represents the evolution of total energy. 
  Dotted and solid lines represent the evolutions of in-slice and out-of-slice 
  components of the kinetic energy. 
  Dot-dashed line represents the evolution of potential energy.}
  \label{fig:time-energy-control}
\end{figure}

Figure \ref{fig:time-energy-control} shows the time evolution of the total energy $E$, 
the kinetic energy $K_{u}$ and $K_{v}$, and the potential energy $P$, which are 
defined by the equations \eqref{kinetic_u} to \eqref{potential}. We observe a loss of
energy after the front is formed, which is due to loss of available potential 
energy caused by numerical dissipation in the $\theta_S$ equation as well as loss of $K_v$ due to numerical dissipation in the $v$ equation. This demonstrates the challenge posed by fronts
in atmosphere models; numerical dissipation is necessary to stabilise advection in the model
but is also altering the large scale circulation.

\begin{figure}[p]
  \centering
  \begin{subfigure}{0.45\hsize}
    \centering
    \includegraphics[width=85mm]{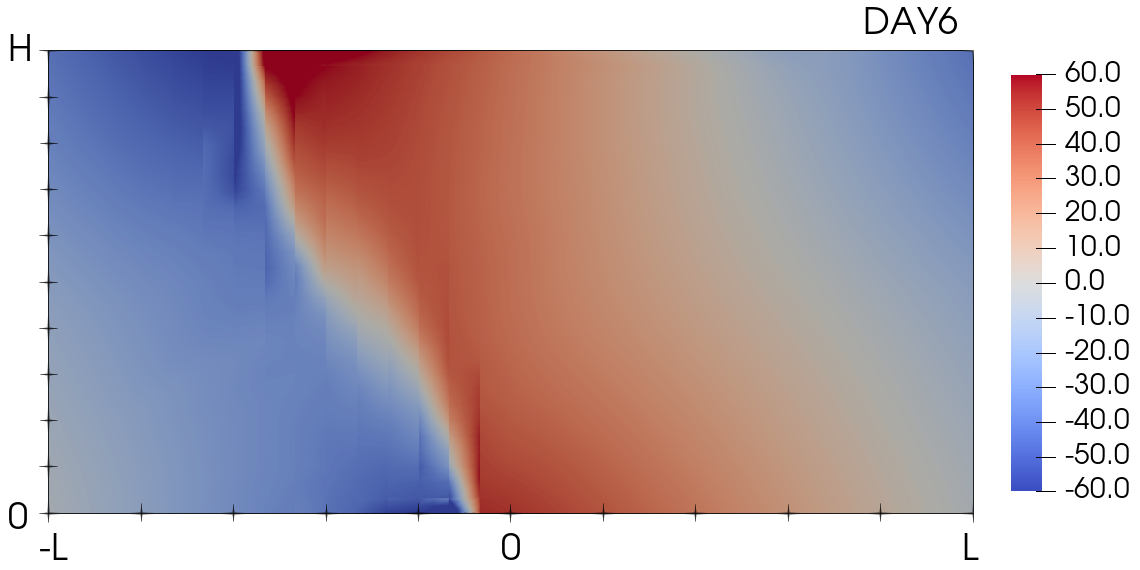}
    \includegraphics[width=85mm]{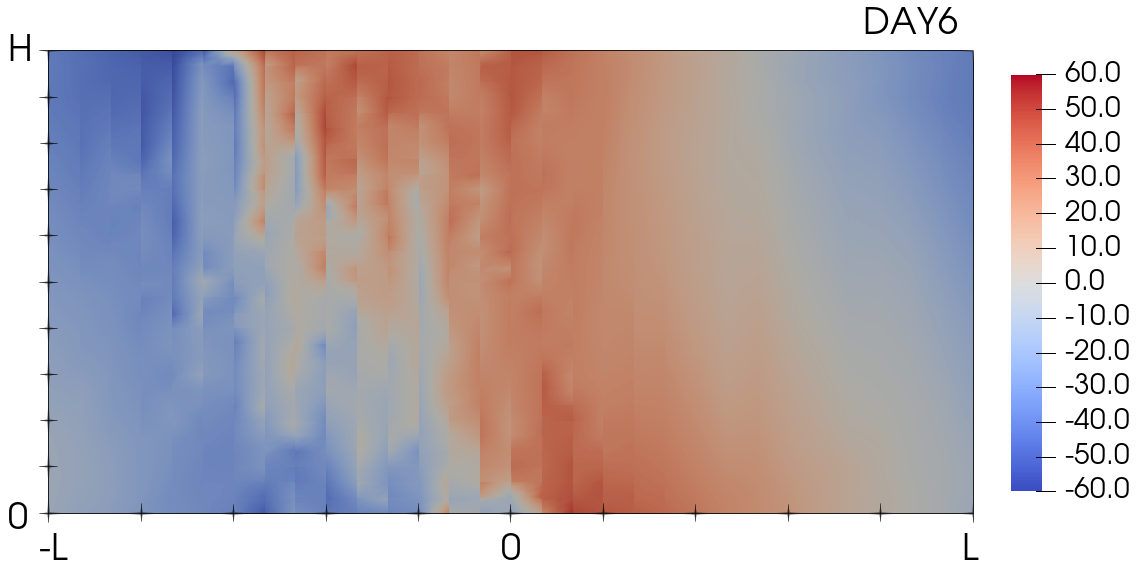}
    \caption{Out-of-slice velocity field}
    \label{fig:velocity-contours-H}
  \end{subfigure}
  \hspace{2em}
  \begin{subfigure}{0.45\hsize}
    \centering
    \includegraphics[width=85mm]{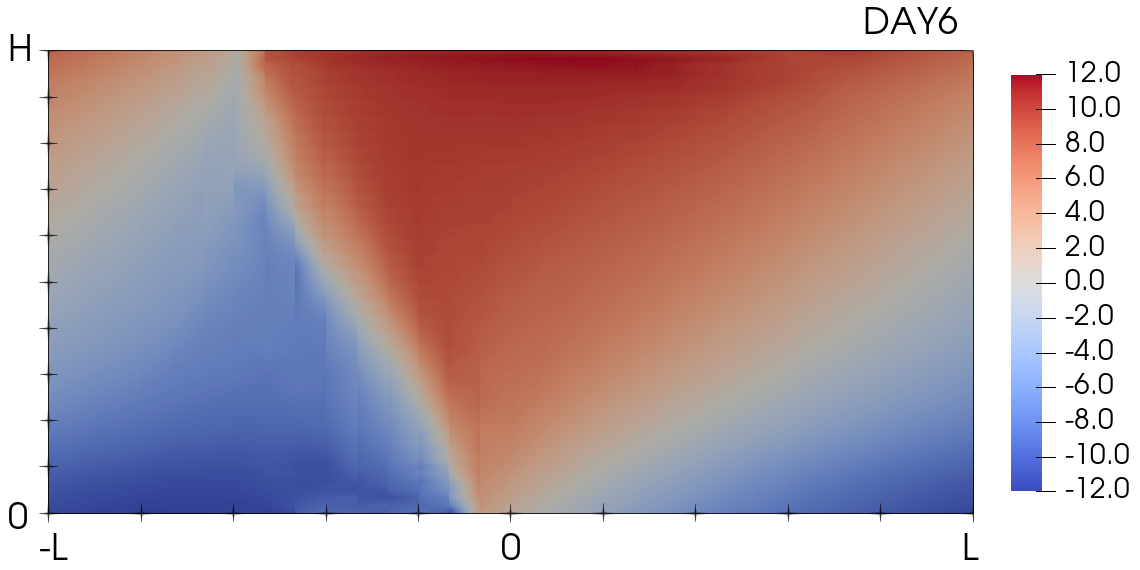}
    \includegraphics[width=85mm]{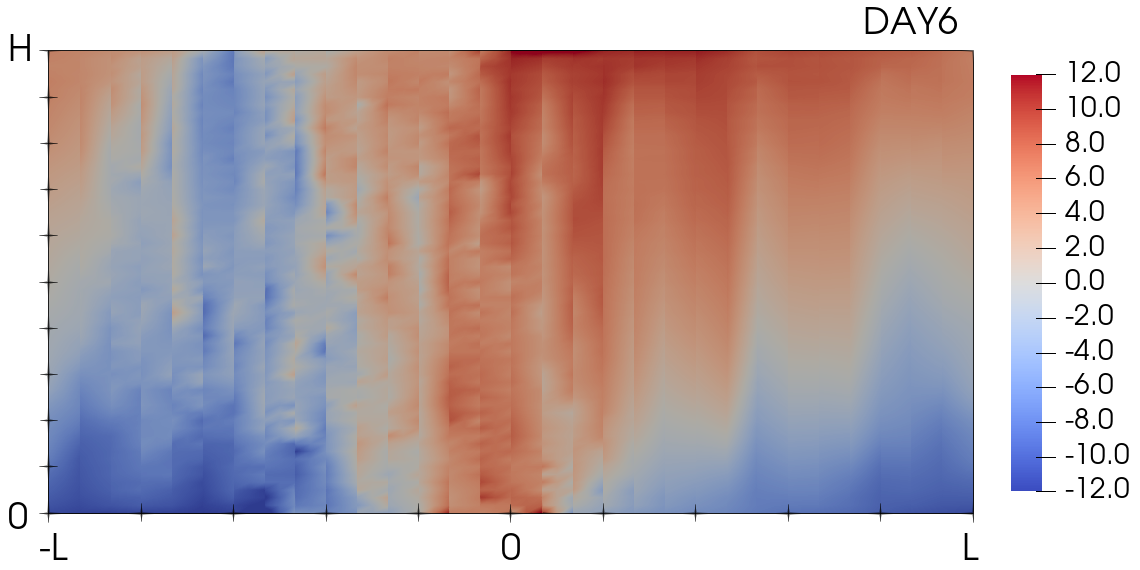}
    \caption{In-slice temperature perturbation field}
    \label{fig:buoyancy-contours-H}
  \end{subfigure}
  \caption{Snapshots of out-of-slice velocity and in-slice temperature perturbation at day 6. top: advective form, bottom: vector invariant form.}
  \label{1fig:slice-contours-H}
\end{figure}

Finally, as an example of how the testcase might be used to compare numerical schemes, we 
replace the advection term $(\MM{u}\cdot\nabla)\MM{u}$ in the velocity equation
with the vector invariant form $(\nabla\times \MM{u})\times \MM{u} + |\MM{u}|^2/2$, using
the upwind discretisation of \citet{cotter2023compatible}. As shown in Figure \ref{1fig:slice-contours-H},
by day 6 the version using the vector invariant form has developed oscillations that
are polluting the solution. We believe that this is related to the ``Hollingsworth instability" 
which is much discussed in the dynamical core literature \citep{hollingsworth1983internal,lazic1984non,gassmann2013global,bell2017numerical, peixoto2018numerical}.
This example requires further investigation, but might provide a useful laboratory to
further understanding how to control that instability.

\section{Conclusion}\label{conclusion}
A new vertical slice test case is developed for compressible nonhydrostatic dynamical cores of atmospheric models, and solutions of a model using an implicit time-stepping method and a compatible finite element method were presented as a reference. The test case can be easily run by adding two terms to a 3D atmosphere model in local area mode, i.e. a Cartesian box, with periodicity and being one cell wide in the $y$-direction. It provides the capability to explore the behaviour of discretisations
in the presence of challenging fronts but in a quasi-2D configuration that allows rapid turnaround
for numerical experiments. The equations for the test case do not correspond to the reduction of a 3D model to a 2D vertical slice configuration, since this is not possible for a compressible model with nonlinear equation of state. However, it serves our purpose since we observe solutions that produce features characteristic of atmospheric fronts, which we would like to challenge our numerical schemes
with.

We demonstrated that the test case produces quasi-periodic lifecycles of fronts which are observed despite the presence of strong discontinuities; these are sustained because the equations have conserved
energy and potential vorticity. Without this formulation, solutions are insufficiently constrained to
drive strong fronts, as was observed by \citet{cullen2008comparison}. The general results of frontogenesis are consistent with results from the incompressible Euler–Boussinesq models, so we have been successful in our goals for the test case. 

We hope that this test case becomes adopted by atmosphere modelling groups, allowing intercomparisons of different discretisations and approaches. One possible use for the testcase is to understand how subgrid schemes, both explicit and implicit, behave in the region of fronts. For example, the compressible version of the $\alpha$-Euler vertical slice model of \citet{cotter2013variational} could be used for this. In addition, the test case could be extended to include the effects of moisture, and boundary layer physics, for example.

To complete the programme of \citet{cullen2007modelling} in this setting, it remains to prove rigorously that the solution of our compressible vertical slice model converges in the semigeostrophic limit to solutions of their semigeostrophic counterpart, following an adaptation of \citet{cullen2003fully} to include the extra terms marked in (\eqref{ueq}-\ref{rhoeq}). This analysis can also be used to obtain a semigeostrophic reference solution, following the approaches of \citet{egan2022new,benamou2024entropic}. We leave this mathematical analysis to future work.

\paragraph{Acknowledgement}
We are grateful for funding from EPSRC via grant EP/R029628/1, and from NERC via grant NE/K012533/1. We would like to thank Thomas Bendall, Mike Cullen, Darryl Holm and Jemma Shipton for their helpful comments on our paper. Any errors are the author's own.

\bibliography{Compressible_eady_front}

\appendix

\section{Variational derivation of the model}

In this appendix we briefly summarise the variational derivation of the vertical slice compressible Eady
model, since it slightly differs from the model presented in Section 5 of \citet{cotter2013variational},
due to the modification to the Exner pressure. The Lagrangian for the system is
\begin{equation}
\ell[\MM{u},v,D,\theta_S] =
\int_\Omega \underbrace{\frac{D}{2}\left(|\MM{u}|^2 + v^2\right) + fDvx}_{\mbox{total kinetic energy}} \qquad -\underbrace{gDz}_{\mbox{gravitational potential energy}} - \underbrace{D\theta_S(c_v\Pi-c_p\Pi_0)}_{\mbox{internal energy}} \diff x\diff z,
\label{eq:lagrangian}
\end{equation}
where $\theta_S$ is the in-slice temperature, transported according to
\begin{equation}
\label{eq:thetaSapp}
\pp{\theta_S}{t} + \MM{u}\cdot \nabla \theta_S + sv = 0,
\end{equation}
with the idea that the total potential temperature is $\theta\approx \theta_S + sy$ in the region
of $y=0$, where $D$ is the in-slice density, transported according to 
\begin{equation}
\label{eq:Dapp}
\pp{D}{t} + \nabla\cdot(\MM{u}D) = 0.
\end{equation}
\eqref{eq:lagrangian} differs from the Lagrangian in \citet{cotter2013variational} by the inclusion of $\Pi_0$, a
constant. Combining the two formulae in \eqref{eq:exner} gives
\begin{equation}
\Pi = \left(\frac{DR\theta_S}{p_0}\right)^{\gamma-1},
\end{equation}
where $\gamma=c_p/c_v$, and we have 
\begin{equation}
\pp{\Pi}{\theta_S} = \frac{\gamma-1}{\theta_S}\Pi = \frac{R}{c_v\theta_S}\Pi, \,
\pp{\Pi}{D} = \frac{R}{c_vD}\Pi.
\end{equation}

The functional derivatives are then
\begin{align}
\dede{\ell}{\MM{u}_S} &= D\MM{u}_S, \\
\dede{\ell}{v} &= D(v + fx), \\
\dede{\ell}{D} &= \frac{1}{2}(\|\MM{u}\|^2 + v^2) + fvx - gz - c_p\theta_S(\Pi - \Pi_0), \\
\dede{\ell}{\theta_S} & = c_pD(\Pi-\Pi_0).
\end{align}
Substitution into Equations 2.9 of \cite{cotter2013variational} then gives
\begin{align}
\nonumber
\pp{D\MM{u}}{t} 
+ \nabla\cdot\left(\MM{u} \otimes D\MM{u}\right) +
(\nabla\MM{u})^T
\cdot D\MM{u}
+ D(v + fx)\nabla v \qquad \qquad \qquad \qquad & \\
\quad - D\nabla\left(\frac{1}{2}(\|\MM{u}\|^2 + v^2) + fvx - gz - c_p\theta_S(\Pi - \Pi_0)\right)
+ c_pD(\Pi-\Pi_0)\nabla\theta_S & = 0, \\
\pp{D(v+fx)}{t} + \nabla\cdot(\MM{u} D(v+fx)) + sc_pD(\Pi-\Pi_0) &= 0,
\end{align}
and further manipulations in combination with Equations \eqref{eq:thetaSapp} and \eqref{eq:Dapp}
lead to Equations (\ref{ueq}-\ref{rhoeq}), noting that $\nabla \Pi_0=0$ (because $\Pi_0$ is a constant).

Appendix B of \cite{cotter2013variational} derives the Lie-Poisson formulation for the vertical slice framework, demonstrating a conserved energy/Hamiltonian given by
\begin{equation}
E = \int_\Omega \dede{\ell}{\MM{u}}\cdot\MM{u}+\dede{\ell}{v}v\diff x\diff z - \ell[\MM{u},v,\theta_S,D].
\end{equation}
For $\ell$ from \eqref{eq:lagrangian}, this becomes
\begin{equation}
E = \int_\Omega {\frac{D}{2}\left(|\MM{u}|^2 + v^2\right)} + {gDz}+ {D\theta_S(c_v\Pi-c_p\Pi_0)} \diff x\diff z,
\end{equation}
which just corresponds to the (relative, i.e. the Coriolis part is removed) kinetic energy
plus the potential energy, as expected. 

Finally, from Equation 2.11 of \cite{cotter2013variational} we have the vertical slice Kelvin circulation theorem,
\begin{equation}
\dd{}{t}\oint_{C(t)}\left(s\frac{1}{D}\dede{\ell}{\MM{u}} - \frac{1}{D}\dede{\ell}{v}\nabla\theta_S\right)\cdot \diff(x,z) = 0,
\end{equation}
where $C(t)$ is a closed loop transported in the $x$-$z$ plane by the in-slice velocity
$\MM{u}$. For $\ell$ from \eqref{eq:lagrangian}, this becomes
\begin{equation}
\dd{}{t}\oint_{C(t)}\left(s\MM{u} - (v + fx)\nabla\theta_S\right)\cdot \diff{(x,z)} = 0.
\end{equation}
After application of Stokes' Theorem and consideration of arbitrary $C(t)$, we obtain
the material conservation of potential vorticity $q$,
\begin{equation}
\pp{q}{t}+\MM{u}\cdot\nabla q = 0, \quad
q = \left(s\curl\MM{u} + \nabla\theta_S \times (\nabla v + f\hat{x})\right)\cdot \hat{y} \diff(x,z),
\end{equation}
where $\hat{x}$ and $\hat{y}$ are the unit vectors pointing in the $x$- and $y$- directions respectively.

\end{document}